# The majority preference relation based on cone preference relations of the decision makers

Alexey O. Zakharov *

St. Petersburg State University, Universitetskii pr., 35, Petrodvorets, St. Petersburg, 198504 Russia

**Abstract**

A multicriteria group choice problem is considered in the paper. The model includes a set of feasible alternatives, a vector criterion, and *n* preference relations of the decision makers (DMs). Each preference relation is a cone relation with corresponding properties. It is considered the majority preference relation, as a cone relation constructed upon the cones of the DMs' preference relations. It is shown how to use and aggregate additional information about the DMs' preference relations in case of three DMs and two criteria. The Pareto set of multicriteria problem with "new" vector criterion forms a group choice, which reduces the Pareto set of initial multicriteria problem.

## 1. Introduction

Problems of choice are considered in deferent areas of human life, and most of them are multicriteria. Wishes and preferences are extremely varied, even opposite to each other. For example, searching the most profitable solution could lead us to the solution with great losses at the same time, and vice versa.

A multicriteria group choice problem includes a group of decision makers (DMs) with individual preferences, which should be aggregated by some method, rule to make a choice[1–8]. The individual preferences can be represented by ranking, estimations of experts, preference relations, utility functions. One of the main problem is how to make consistent individual preferences, what should be considered as a "good" decision. A. Bergson, P. Samuelson proposed a social welfare function as a model of individual aggregation, K. J. Arrow formulated 5 axioms for individual preference relations, which restricted the "reasonable" group choice, and Arrow paradox is said that there is not an aggregation rule satisfied this axioms[1]. M.A. Aizerman and F.T. Aleskerov proposed the axioms of "reasonable" group (collective) choice for individual choice functions[2]. B.G. Mirkin investigated problems of preference descriptions, consistency principle of individual preferences, properties of majority rule and its modifications[4].

A model considered in the paper includes a set of feasible alternatives, a vector criterion reflecting the goals of the group, and preference relations of DMs. Information about preference relation of the DM (the individual preferences) is given by "quanta" of information[9, 10], so preference relation of the DM is characterized by a convex pointed cone, which contains nonnegative orthant, and does not contain the origin. For aggregating the individual preference relations is used the majority preference relation, which takes into account at least half intensions of members of the group. It is proved that this relation is characterized by the unions of the intersections of different subsets selected from the set of cones associated with the DMs' preference relations. Generally, derived in this way cone can be either convex or not. The property of convex is equivalent to transitivity of corresponding relation. Thus, the goal is to construct a convex part of the majority preference relation cone, and it gives the transitive part of this relation.

---

\* Corresponding author. Tel.: +7-812-428-48-68.
  *E-mail address:* a.zakharov@spbu.ru

When there is not any "quantum" of information, the Pareto set forms the group choice, i. e. each DM has not any additional information. For three DMs and two components of vector criterion it is considered two cases: when each DM have only one "quantum" of information, and when it has two "quanta" of information. And it is established what consists the group choice in each case: it is the Pareto set of multicriteria problem with "new" vector criteria, and it belongs to the Pareto set of initial multicriteria problem. Note that in one considered situation this choice is not unique. It happens, when the cone of the majority preference relation is not convex.

## 2. Multicriteria choice model of the DM

Consider the group consists of $n$ DMs: $DM_1, \ldots, DM_n$. A set of feasible alternatives, solutions $X \subseteq R^k$ is a set of variants among which the choice should be made by the group. A vector criterion $\mathbf{f}$, defined on $X$, $Y = \mathbf{f}(X)$, reflects goals, intensions of the DMs, and it is the same for each DM. Now let us formulate the main objects of multicriteria group choice problem:
- a set of feasible vectors $Y$;
- $n$ preference relations $\succ_1^Y, \ldots, \succ_n^Y$, defined on $Y$, of $DM_1, \ldots, DM_n$.

This relations $\succ_1^Y, \ldots, \succ_n^Y$ also reflects intensions all DMs but they are not equal to each other, i. e. for any $l \in \{1, \ldots, n\}$ $DM_l$ has its own preference relation $\succ_l^Y$, defined on $Y$. The expression $\mathbf{y}^{(1)} \succ_l^Y \mathbf{y}^{(2)}$ means that when considering two vectors $\mathbf{y}^{(1)}, \mathbf{y}^{(2)}$ the $DM_l$ chooses $\mathbf{y}^{(1)}$ and rejects $\mathbf{y}^{(2)}$.

Suppose that behavior of any $DM_l$ is restricted by the axioms of "reasonable" individual choice[9]. They are as follows. It is assumed that when considering two possible variants the excluded vector could not be selected from the whole set $Y$ as well. To have the opportunity to compare the random vectors (not only from the set $Y$), and thus provide information about the preference relation, it is introduced an irreflexive and transitive continuation $\succ_l$ to the entire criterion space $R^m$ for the latter. In addition, all components $f_1, \ldots, f_m$ of vector criterion are compatible with preference relation $\succ_l$. The $i$th criterion $f_i$ is called compatible with preference relation $\succ_l$ if for any vectors $\mathbf{y}^{(1)}, \mathbf{y}^{(2)}$ of the space $R^m$ such that $\mathbf{y}^{(1)} = (y_1^{(1)}, \ldots, y_{i-1}^{(1)}, y_i^{(1)}, y_{i+1}^{(1)}, \ldots, y_m^{(1)})$, $\mathbf{y}^{(2)} = (y_1^{(1)}, \ldots, y_{i-1}^{(1)}, y_i^{(2)}, y_{i+1}^{(1)}, \ldots, y_m^{(1)})$, $y_i^{(1)} > y_i^{(2)}$, the relation $\mathbf{y}^{(1)} \succ_l \mathbf{y}^{(2)}$ holds. It means that the DM is interested in increasing value of each criterion while values of others criteria are constant. Thus, it is outlined the class of multicriteria individual choice problem, and it is developed the Pareto set reduction axiomatic approach by V.D. Noghin applied to this class[9–12].

It is shown[9] that under assumption the axioms of "reasonable" choice being valid, the Edgeworth – Pareto principle for any $DM_l$ is hold. It is said that the "best" solution should be chosen only within the set of Pareto-optimal vectors (the Pareto set) $P(Y)$. Here, $P(Y) = \{\mathbf{y}^* \in Y \mid \bar{\exists} \mathbf{y} \in Y : \mathbf{y} \geq \mathbf{y}^*\}$, also one can consider the set of Pareto-optimal solutions $P_f(X) = \{\mathbf{x}^* \in X \mid \bar{\exists} \mathbf{x} \in X : \mathbf{f}(\mathbf{x}) \geq \mathbf{f}(\mathbf{x}^*)\}$. The relation $\mathbf{y} \geq \mathbf{y}^*$ means that the relations $y_i \geq y_i^*$ are valid for all $i = 1, ..., m$, and $\mathbf{y} \neq \mathbf{y}^*$. It is usually called the Pareto relation.

To specify the individual choice within the Pareto set it is given some additional information about the preference relation $\succ_l$ of the $DM_l$. It is introduced notation of "quantum" of information[9, 11, 12]. Following it let us recall the definition.

**Definition 1.** It is said that we have a "quantum" of information about the $DM_l$'s preference relation $\succ_l$ with groups of criteria $A$ and $B$ and with two sets of positive parameters $w_i^{(l)}$ for all $i \in A$ and $w_j^{(l)}$ for all $j \in B$ if for all vectors $\mathbf{y}^{(1)}, \mathbf{y}^{(2)} \in R^m$ such that

$$y_i^{(1)} - y_i^{(2)} = w_i^{(l)} > 0 \quad \forall i \in A, \quad y_j^{(2)} - y_j^{(1)} = w_j^{(l)} > 0 \quad \forall j \in B, \quad y_s^{(1)} = y_s^{(2)} \quad \forall s \in I \setminus (A \cup B),$$

where $I = \{1, \ldots, m\}$, $A, B \subset I$, $A \neq \emptyset$, $B \neq \emptyset$, $A \cap B = \emptyset$, the following relation is valid: $\mathbf{y}^{(1)} \succ_l \mathbf{y}^{(2)}$. In such case the group of criteria $A$ is called more important, and the group $B$ is called less important with given positive parameters $w_i^{(l)}$ (profits), $w_j^{(l)}$ (losses).

Add one more condition to the list of axioms of "reasonable" choice: the preference relation $\succ_l$ for any $l \in \{1, \ldots, n\}$ is invariant under a linear positive transformation[9]. Taking it into account it is shown[8] that in Definition 1 the vectors $\mathbf{y}^{(1)}$, $\mathbf{y}^{(2)}$ can be assumed to be fixed, in particular, we can put $y_i^{(1)} = w_i^{(l)}$ $\forall i \in A$, $y_j^{(1)} = -w_j^{(l)}$ $\forall j \in A$, $y_s^{(1)} = 0$ $\forall s \in I \setminus (A \cup B)$, $\mathbf{y}^{(2)} = \mathbf{0}_m$. So, the existence of "quantum" of information means that it is given a vector $\mathbf{y} \in N^m$ ($N^m = R^m \setminus (R_+^m \cup (-R_+^m) \cup \{0_m\})$) such that $\mathbf{y} \succ_l \mathbf{0}_m$.

In papers[9–12] in order to specify the individual choice it is shown how to use collections of "quantum" of information of different types, which according to the aforementioned could be defined by the sequence of vectors $\mathbf{y}^{(s)} \in N^m$, $s = 1, \ldots, p_l$. And if there exists the preference relation $\succ_l$ satisfied the axioms of "reasonable" choice such that $\mathbf{y}^{(s)} \succ_l \mathbf{0}_m$ for any $l \in \{1, \ldots, n\}$, then such information about the $DM_l$'s preference relation is called consistent[9]. The inconsistent information could not be used in decision making process. The criteria of consistency are derived[9]. The main idea of the Pareto set reduction axiomatic approach is to construct a set, which will restrict a "new" bound of individual choice using the given information. And such bound is more precise, than the Pareto set.

## 3. The majority preference relation

**Definition 2.** A binary relation $\Re$ defined on $R^m$ is called a cone relation if it is existed a cone $K$ such that the following equivalent is hold for any $\mathbf{y}^{(1)}, \mathbf{y}^{(2)} \in R^m$: $\mathbf{y}^{(1)} \Re \mathbf{y}^{(2)} \Leftrightarrow \mathbf{y}^{(1)} - \mathbf{y}^{(2)} \in K$. Remark that the inclusion $\mathbf{y}^{(1)} - \mathbf{y}^{(2)} \in K$ is the same to $\mathbf{y}^{(1)} \in \mathbf{y}^{(2)} + K$.

It is proved[9] that the preference relation which satisfies the axioms of "reasonable" choice is a cone relation with a pointed convex cone (without the origin $\mathbf{0}_m$) that contains the nonnegative orthant $R_+^m = \{\mathbf{y} \in R^m \mid \mathbf{y} \geq \mathbf{0}_m\}$. For that reason the preference relation $\succ_l$ of the $DM_l$ for any $l = 1, \ldots, n$ is a cone relation with abovementioned preferences. Denote such cone $K_l$ for each $\succ_l$.

Consider a group preference relation, which takes into account the intentions of at least a half of the group. For example, if for some vectors $\mathbf{y}^{(1)}, \mathbf{y}^{(2)} \in R^m$ the relations $\mathbf{y}^{(1)} \succ_1 \mathbf{y}^{(2)}$, $\mathbf{y}^{(1)} \succ_2 \mathbf{y}^{(2)}$, $\ldots$, $\mathbf{y}^{(1)} \succ_q \mathbf{y}^{(2)}$ are valid, where $q$ is greater than $n/2$, then the relation $\mathbf{y}^{(1)} \succ \mathbf{y}^{(2)}$ holds.

**Definition 3.** Let us call the group preference relation $\succ$ defined on $R^m$ the *majority preference relation* if for any vectors $\mathbf{y}^{(1)}, \mathbf{y}^{(2)} \in R^m$ the relation $\mathbf{y}^{(1)} \succ \mathbf{y}^{(2)}$ is equivalent to the existence of such subset $\{l_1, \ldots, l_p\} \subset \{1, \ldots, n\}$ that the relations $\mathbf{y}^{(1)} \succ_{l_j} \mathbf{y}^{(2)}$ are valid for all $j = 1, \ldots, p$, where $p = n/2$, if $p$ is even, and $p = \left[\frac{n}{2}\right] + 1$, if $p$ is odd.

Here for any number $a$ by $[a]$ we denote an integer part of number $a$. Consider the cone

$$K = \bigcup_{l=1}^{C_n^p} \bigcap_{j=1}^{p} K_{lj} \qquad (1)$$

where $\{K_{l_1}, \ldots, K_{l_p}\}$ is the subset of the set of the cones $\{K_1, \ldots, K_n\}$ of preference relations $\succ_l$, $l = 1, \ldots, n$. Remark that for any $l, k \in \{1, \ldots, n\}$, $l \neq k$, subsets $\{K_{l_1}, \ldots, K_{l_p}\}$ and $\{K_{k_1}, \ldots, K_{k_p}\}$ do not coincide, $C_n^p$ – $p$-combination of a set $\{1, \ldots, n\}$. Obviously, $K$ is a cone as unions and intersections of cones.

**Lemma 1.** *The majority preference relation $\succ$ is a cone relation with cone $K$.*

**Proof.** It means that it is sufficient to establish the following equivalent: $\mathbf{y}^{(1)} \succ \mathbf{y}^{(2)} \Leftrightarrow \mathbf{y}^{(1)} - \mathbf{y}^{(2)} \in K$.

Necessity. According to Definition 3 if the relation $\mathbf{y}^{(1)} \succ \mathbf{y}^{(2)}$ holds, then there exists such subset $\{l_1, \ldots, l_p\} \subset \{1, \ldots, n\}$ that the relations $\mathbf{y}^{(1)} \succ_{l_j} \mathbf{y}^{(2)}$ are valid for all $j = 1, \ldots, p$, where $p = n/2$, if $p$ is even, and $p = \left[\frac{n}{2}\right] + 1$, if $p$ is odd. Since the preference relation $\succ_s$ of the $DM_s$ for any $s = 1, \ldots, n$ is a cone relation with the appropriate cone $K_s$, there exists a subset $\{K_{l_1}, \ldots, K_{l_p}\} \subset \{K_1, \ldots, K_n\}$ that $\mathbf{y}^{(1)} - \mathbf{y}^{(2)} \in K_{l_j}$ for all $j = 1, \ldots, p$. This implies the inclusion $\mathbf{y}^{(1)} - \mathbf{y}^{(2)} \in \bigcap_{j=1}^{p} K_{l_j}$, and we have $\mathbf{y}^{(1)} - \mathbf{y}^{(2)} \in K$.

**Sufficiency.** Let the inclusion $\mathbf{y}^{(1)} - \mathbf{y}^{(2)} \in K$ is valid. It means that there exists such number $l \in \{1, \ldots, C_n^p\}$ that $\mathbf{y}^{(1)} - \mathbf{y}^{(2)} \in \bigcap_{j=1}^{p} K_{lj}$, so we obtain that the inclusion $\mathbf{y}^{(1)} - \mathbf{y}^{(2)} \in K_{lj}$ is true for any $j \in \{1, \ldots, p\}$. This implies the relation $\mathbf{y}^{(1)} \succ_{l_j} \mathbf{y}^{(2)}$, where $\succ_{l_j}$ is the relation appropriate to the cone $K_{lj}$. As a result we have that for some number $l \in \{1, \ldots, C_n^p\}$ we specified such subset $\{l_1, \ldots, l_p\} \subset \{1, \ldots, n\}$ that the relations $\mathbf{y}^{(1)} \succ_{l_j} \mathbf{y}^{(2)}$ are valid for any $j \in \{1, \ldots, p\}$. According to Definition 3 the relation $\mathbf{y}^{(1)} \succ \mathbf{y}^{(2)}$ holds. Lemma 1 is proved.

We obtained that the group preference relation defined in Definition 3 is a cone relation. Now consider the properties of this relation and its cone.

**Lemma 2.** *The majority preference relation $\succ$ is an irreflexive relation, invariant under a linear positive transformation, and its cone K contains the nonnegative orthant $R_+^m$, and does not contain the origin $\mathbf{0}_m$.*

**Proof.** a) Let us prove that the relation $\succ$ is invariant under a linear positive transformation. Consider two arbitrary vectors $\mathbf{y}^{(1)}, \mathbf{y}^{(2)} \in R^m$ such that $\mathbf{y}^{(1)} \succ \mathbf{y}^{(2)}$ holds. Since the relation $\succ$ is a cone relation with cone $K$ defined in (1), we have $\mathbf{y}^{(1)} - \mathbf{y}^{(2)} \in K$. This is equivalent to the inclusion $\alpha \mathbf{y}^{(1)} - \alpha \mathbf{y}^{(2)} \in K$ for any nonzero $\alpha$, because the set $K$ is a cone. For any arbitrary vector $\mathbf{c} \in R^m$ we obtain $\alpha \mathbf{y}^{(1)} - \alpha \mathbf{y}^{(2)} \in K \Leftrightarrow (\alpha \mathbf{y}^{(1)} + \mathbf{c}) - (\alpha \mathbf{y}^{(2)} + \mathbf{c}) \in K$. The right hand side inclusion is equivalent to the relation $(\alpha \mathbf{y}^{(1)} + \mathbf{c}) \succ (\alpha \mathbf{y}^{(2)} + \mathbf{c})$. As a result, the relation $\succ$ is invariant under a linear positive transformation.

b) It is easy to see from Definition 3 that the relation $\succ$ is irreflexive due to the property of irrefexivity of all relations $\succ_1, \ldots, \succ_n$ (it can be proved by contradiction).

c) Due to the inclusions $R_+^m \subset K_l$ and $0_m \notin K_l$ for any $l \in \{1, \ldots, n\}$, it is obvious that $R_+^m \subseteq K$ and $0_m \notin K$. Lemma 2 is proved.

In spite of cones $K_l$ for all $l \in \{1, \ldots, n\}$ are pointed and convex, in general, the cone $K$ is not pointed and convex. The property of convex for cone means that corresponding cone relation is transitive[9]. One can obtain it using the definition of the transitivity and the fact that for any convex cone $\hat{K}$ and for any arbitrary vectors $\mathbf{y}^{(1)}, \mathbf{y}^{(2)} \in \hat{K}$ the inclusion $\mathbf{y}^{(1)} + \mathbf{y}^{(2)} \in \hat{K}$ is hold. Let $\hat{K}$ be the convex cone such that $\hat{K} \subseteq K$, call this cone $\hat{K}$ a *transitive part* of the cone $K$. And if it does not exist a convex cone $\tilde{K} \subseteq K$ that $\tilde{K} \subset \hat{K}$, $\tilde{K} \neq \hat{K}$, then call this cone $\hat{K}$ a *maximum transitive part* of the cone $K$. In general, it is not a unique for particular cone $K$. Similarly, let us call the corresponding cone relations $\succ_{tr}$ and $\succ_{\max tr}$ a *transitive part* and a *maximum transitive part* of the majority preference relation $\succ$.

It is obvious that $R_+^m \subseteq \hat{K}$, $0_m \notin \hat{K}$, then the relation $\succ_{tr}$ is irreflexive, transitive, and invariant under a linear positive transformation. Denote by $\text{Ndom}_{\succ_{tr}}(Y)$ the set of nondominated vectors of the set $Y$ according to an arbitrary transitive part $\succ_{tr}$ of the relation $\succ$, i. e. $\text{Ndom}_{\succ_{tr}}(Y) = \{\mathbf{y}^* \in Y \mid \bar{\exists} \mathbf{y} \in Y : \mathbf{y} - \mathbf{y}^* \in \hat{K}\}$. The inclusion $R_+^m \subseteq \hat{K}$ implies the inclusion $\text{Ndom}_{\succ_{tr}}(Y) \subseteq P(Y)$.

If for any $l \in \{1, \ldots, n\}$ $DM_l$ has not any additional information about its preference relation $\succ_l$, i. e. there is not "quantum" of information about the relation $\succ_l$, then the corresponding cone $K_l$ coincides with the nonnegative orthant $R_+^m$. And we have that the cone $K = R_+^m$. In such case all three relations $\succ$, $\succ_{tr}$, and $\geq$ coincide, and all vectors of the Pareto set form a group choice. But the reducing this choice is arisen in practice. The goal is to construct a transitive part $\hat{K}$ of the cone $K$ using "quanta" of information about the DMs' preference relations. The set of nondominated vectors $\text{Ndom}_{\succ_{tr}}(Y)$ of the set $Y$ according to this transitive part $\hat{K}$ will be considered as a group choice based on the given "quanta" of information, and $\text{Ndom}_{\succ_{tr}}(Y) \subseteq P(Y)$. So it reduces the bounds of group choice.

## 4. Group choice using the majority preference relation in case of 3 DMs and m = 2

*4.1. Case of one "quantum" of information for each DM*

Consider the group consists of 3 DMs, each of them is associated to its preference relation with corresponding cone: $DM_1$ with $\succ_1$ (cone $K_1$), $DM_2$ with $\succ_2$ (cone $K_2$), $DM_3$ with $\succ_3$ (cone $K_3$). Assume that preference relations $\succ_1$, $\succ_2$, and $\succ_3$ are satisfied the axioms of "reasonable" choice, so cone $K_l$ is convex and pointed, $R_+^2 \subset K_l$, and $\mathbf{0}_2 \notin K_l$ for any $l \in \{1, 2, 3\}$. Let $\mathbf{f} = (f_1, f_2)$, it means that $m = 2$.

According to Definition 3 and Lemma 1 the majority preference relation $\succ$ is a cone relation with a cone
$$K = (K_1 \cap K_2) \cup (K_1 \cap K_3) \cup (K_2 \cap K_3). \tag{2}$$

**Lemma 3.** *If the cone $K$ defined in (2) is convex, then $\hat{K} = K$. If the cone $K$ is not convex, then an arbitrary transitive part is a cone $\hat{K} = \text{cone}\{\mathbf{e}^1, \mathbf{e}^2, \mathbf{u}^{(1)}, \mathbf{u}^{(2)}\} \setminus \{\mathbf{0}_2\}$, where the vectors $\mathbf{u}^{(1)} = (u_1^{(1)}, -u_2^{(1)})^T$, $\mathbf{u}^{(2)} = (-u_1^{(2)}, u_2^{(2)})^T$, $\mathbf{u}^{(1)}, \mathbf{u}^{(2)} \in K$, with components $u_1^{(1)}, u_2^{(1)}, u_1^{(2)}, u_2^{(2)} > 0$ satisfy the inequality*
$$u_1^{(1)} u_2^{(2)} - u_2^{(1)} u_1^{(2)} > 0. \tag{3}$$
*And it does not exist the maximum transitive part of the cone $K$.*

**Proof.** In case the cone $K$ is convex, then, obviously, $\hat{K} = K$. Consider the case, when the cone $K$ is not convex. The cone $K$ does not contain the set $(-R_+^2)$. If we assume that any arbitrary vector $\mathbf{y} \in (-R_+^2)$ belongs to the cone $K$, then there exists such number $i \in \{1, 2, 3\}$ that for cone $K_i$ the following inclusion is hold: $\mathbf{y} \in K_i$. It implies the condition $\mathbf{y} \succ_i \mathbf{0}_m$, which contradicts with the axioms of "reasonable" choice of the $DM_i$. Thus, $(-R_+^2) \not\subset K$.

According to Lemma 2 the inclusion $R_+^m \subseteq K$ is hold. If $K = R_+^m$, then, obviously, the cone $K$ is convex. It should be noted that the case $K = R_+^m$ has not practical sense, because this implies that the relation $\succ$ coincides with the Pareto relation $\geq$.

As a result, due to $(-R_+^2) \not\subset K$ and $R_+^m \subseteq K$, there exist two vectors $\mathbf{u}^{(1)} = (u_1^{(1)}, -u_2^{(1)})^T$, $\mathbf{u}^{(2)} = (-u_1^{(2)}, u_2^{(2)})^T$ with components $u_1^{(1)}, u_2^{(1)}, u_1^{(2)}, u_2^{(2)} > 0$, which belong to $K$. Consider the cone $\hat{K} = \text{cone}\{\mathbf{e}^1, \mathbf{e}^2, \mathbf{u}^{(1)}\} \cup \text{cone}\{\mathbf{e}^2, \mathbf{u}^{(2)}\}$ without the origin $\mathbf{0}_2$.

Let us prove that the cone $\hat{K}$ is convex if and only if (3). Consider the line $l$ passing through the vector $\mathbf{u}^{(2)}$, its equation is $\langle \mathbf{n}, \mathbf{y} \rangle = 0$, where $\mathbf{n} = (u_2^{(2)}, u_1^{(2)})^T$ is a normal vector. Note that in case $\mathbf{u}^{(1)} \in l$ the cone $\hat{K}$ is not convex due to $\mathbf{0}_2 \notin \hat{K}$, $\mathbf{0}_2 \in l$. If the vector $\mathbf{u}^{(1)}$ generates an acute angle with the normal vector $\mathbf{n}$, i. e. $\langle \mathbf{n}, \mathbf{u}^{(1)} \rangle = u_1^{(1)} u_2^{(2)} - u_2^{(1)} u_1^{(2)} > 0$, then the cone $\hat{K}$ is convex, and it is generated by the vectors $\mathbf{e}^1, \mathbf{e}^2, \mathbf{u}^{(1)}, \mathbf{u}^{(2)}$, and does not contain the origin $\mathbf{0}_2$. We can conclude that the cone $\hat{K}$ in case (3) is an arbitrary transitive part of the cone $K$.

As assumed above the vectors $\mathbf{u}^{(1)}, \mathbf{u}^{(2)} \in K$, and the cone $K$ is not convex. Hence, there be one of two possibilities: 1) there exists such vector $\mathbf{v}^{(1)} = (v_1^{(1)}, -v_2^{(1)}) \in K$, $v_1^{(1)}, v_2^{(1)} > 0$, that $v_1^{(1)} u_2^{(2)} - v_2^{(1)} u_1^{(2)} \leq 0$, or 2) there exists such vector $\mathbf{v}^{(2)} = (-v_1^{(2)}, v_2^{(2)}) \in K$, $v_1^{(2)}, v_2^{(2)} > 0$, that $u_1^{(1)} v_2^{(2)} - u_2^{(1)} v_1^{(2)} \leq 0$. Consider the second one, for the first one arguments are similar. Obviously, $(\text{cone}\{\mathbf{e}^1, \mathbf{e}^2, \mathbf{u}^{(1)}\} \cup \text{cone}\{\mathbf{e}^2, \mathbf{v}^{(2)}\}) \setminus \{\mathbf{0}_2\} \subseteq K$.

There exists such number $\varepsilon > 0$ that
$$u_2^{(2)} > u_2^{(2)} - \varepsilon > \frac{u_1^{(2)} u_2^{(1)}}{u_1^{(1)}}. \tag{4}$$

This implies the inequality $u_1^{(1)}(u_2^{(2)} - \varepsilon) - u_2^{(1)} u_1^{(2)} > 0$. So, the cone $\hat{K}_\varepsilon = \text{cone}\{\mathbf{e}^1, \mathbf{e}^2, \mathbf{u}^{(1)}\} \cup \text{cone}\{\mathbf{e}^2, \mathbf{u}^{(\varepsilon)}\}$ without the origin $\mathbf{0}_2$, where $\mathbf{u}^{(\varepsilon)} = (-u_1^{(2)}, u_2^{(2)} - \varepsilon)^T$, is convex. It is easy to check that the equality $\mathbf{u}^{(2)} = \mathbf{u}^{(\varepsilon)} + \varepsilon \mathbf{e}^2$ is true. It means that the inclusion $\hat{K} \subset \hat{K}_\varepsilon$ is hold. Using (4) and $u_1^{(1)} v_2^{(2)} - u_2^{(1)} v_1^{(2)} \leq 0$ one can check that the linear combination $\mathbf{u}^{(\varepsilon)} = \alpha_1 \mathbf{v}^{(2)} + \alpha_2 \mathbf{e}^2$ is valid for some positive numbers $\alpha_1$, $\alpha_2$. Thus, $\mathbf{u}^{(\varepsilon)} \in \text{cone}\{\mathbf{e}^2, \mathbf{v}^{(2)}\}$. And we obtain the inclusion $\hat{K}_\varepsilon \subset K$.

We have shown that for any arbitrary transitive part $\hat{K}$ there exists a convex cone $\hat{K}_\varepsilon \subset K$ that $\hat{K} \subset \hat{K}_\varepsilon$, $\hat{K}_\varepsilon \neq \hat{K}$. It means that if the cone $K$ is not convex, then it does not exist the maximum transitive part of it. Lemma 3 is proved.

Now, let each $DM_l$ give a "quantum" of information about its preference relation $\succ_l$, $l = 1, 2, 3$. Due to this, we have 3 vectors $\mathbf{y}^{(1)}, \mathbf{y}^{(2)}, \mathbf{y}^{(3)} \in N^2$ such that $\mathbf{y}^{(1)} \succ_1 \mathbf{0}_2$, $\mathbf{y}^{(2)} \succ_2 \mathbf{0}_2$, and $\mathbf{y}^{(3)} \succ_3 \mathbf{0}_2$. There are two and only two situations: (I) the index of positive component is the same for all vectors, for example $y_1^{(l)} > 0$ for any $l \in \{1, 2, 3\}$; (II) only two vectors have positive component with the same index, for example $y_1^{(1)} > 0$, $y_1^{(2)} > 0$, $y_2^{(3)} > 0$. Firstly, we consider the situation (I). Denote $W(l, s) = w_1^{(l)} w_2^{(s)} - w_2^{(l)} w_1^{(s)}$.

**Theorem 1.** *Let the situation (I) is hold, and the vectors*

$$\mathbf{y}^{(1)} = \begin{pmatrix} w_1^{(1)} \\ -w_2^{(1)} \end{pmatrix}, \ \mathbf{y}^{(2)} = \begin{pmatrix} w_1^{(2)} \\ -w_2^{(2)} \end{pmatrix}, \ and \ \mathbf{y}^{(3)} = \begin{pmatrix} w_1^{(3)} \\ -w_2^{(3)} \end{pmatrix}$$

*such that $\mathbf{y}^{(1)} \succ_1 \mathbf{0}_2$, $\mathbf{y}^{(2)} \succ_2 \mathbf{0}_2$, and $\mathbf{y}^{(3)} \succ_3 \mathbf{0}_2$ are given. It is assumed that these vectors are not codirectional. Then $\hat{P}_\mathbf{g}(Y) \subseteq P(Y)$, where $\hat{P}_\mathbf{g}(Y) = f(P_\mathbf{g}(X))$, and $g_1 = f_1$, $g_2 = w_2^{(s)} f_1 + w_1^{(s)} f_2$. Here the index $s \in \{1, 2, 3\}$ is selected as follows: it satisfies the inequalities $W(l, s) > 0$, $W(s, k) > 0$ for some indices $l, k \in \{1, 2, 3\}$, $l \neq k$, $l \neq s$, $k \neq s$. Note that there is only one such index $s$ for particular components of vectors $\mathbf{y}^{(1)}$, $\mathbf{y}^{(2)}$, and $\mathbf{y}^{(3)}$.*

**Proof.** Let $M_l = \text{cone}\{\mathbf{e}^1, \mathbf{e}^2, \mathbf{y}^{(l)}\} \setminus \{\mathbf{0}_2\}$ for any $l \in \{1, 2, 3\}$.

Let us show that for cones $M_l$, $M_s$ for any $l, s \in \{1, 2, 3\}$, $l \neq s$, it is valid either inclusion $M_l \subset M_s$ or inclusion $M_s \subset M_l$. Suppose that the condition $W(l, s) > 0$ is hold. Consider an arbitrary vector $\mathbf{y} \in M_l$. According to the definition of cone $M_l$ we derive that it is existed nonnegative numbers $\lambda_1$, $\lambda_2$, and $\mu_l$ such that $\mathbf{y} = \lambda_1 \mathbf{e}^1 + \lambda_2 \mathbf{e}^2 + \mu_l \mathbf{y}^{(l)}$. It is easy to check that the vector $\mathbf{y}^{(l)}$ can be represented as a linear combination $\mathbf{y}^{(l)} = \alpha_1 \mathbf{e}^1 + \alpha_2 \mathbf{y}^{(s)}$, where $\alpha_1 = W(l, s) / w_2^{(s)}$, $\alpha_2 = w_2^{(l)} / w_2^{(s)}$. In other words, we establish that $\mathbf{y} = \bar{\lambda}_1 \mathbf{e}^1 + \lambda_2 \mathbf{e}^2 + \bar{\mu}_l \mathbf{y}^{(s)}$, where due to $W(l, s) > 0$ numbers $\bar{\lambda}_1$, $\lambda_2$, and $\bar{\mu}_l$ are nonnegative. It means that $\mathbf{y} \in M_s$, so we have $M_l \subset M_s$ in case of $W(l, s) > 0$. Similarly, it can be obtained that the inclusion $M_s \subset M_l$ is valid, if $W(l, s) < 0$. It should be noted that here we consider only strong inclusions $M_l \subset M_s$ and $M_l \subset M_s$. The equality $M_l = M_s$ is possible if and only if $W(l, s) = 0$.

Summarizing the previous analysis we can make a conclusion. The cones $M_1$, $M_2$, and $M_3$ satisfy the inclusions $M_l \subset M_s \subset M_k$, where indices $l, s, k \in \{1, 2, 3\}$, $l \neq s \neq k$. Obviously, in such case the inequalities $W(l, s) > 0$ and $W(s, k) > 0$ are hold. (Remark that $W(l, s) > 0$ and $W(s, k) > 0$ leads to $W(l, k) > 0$). Due to this, the intersections of the cones are as follows: $M_l \cap M_s = M_l$, $M_l \cap M_k = M_l$, $M_s \cap M_k = M_s$. As a result we obtain $M = M_s$. It means that the cone $M$ in situation (I) is convex, also it contains nonnegative orthant $R^2_+$, and does not contain the origin $\mathbf{0}_2$. One can conclude that the cone relation $\succ_{M_s}$ with the cone $M_s$ is irreflexive, transitive, and invariant under a linear positive transformation[9].

As a result, the relation $\succ_{M_s}$ is a transitive part of the majority preference relation $\succ$. From the inclusion $R^2_+ \subset M_s$ we have $\text{Ndom}_{M_s}(Y) \subseteq P(Y)$, where $\text{Ndom}_{M_s}(Y)$ − is the set of nondominated vectors of the set $Y$ according to the relation $\succ_{M_s}$: $\text{Ndom}_{M_s}(Y) = \{\mathbf{y}^* \in Y \mid \bar{\exists} \mathbf{y} \in Y: \mathbf{y} - \mathbf{y}^* \in M_s\}$.

One can prove that the cone $M_s$ coincides with the set of all nonzero solutions of the following system of linear inequalities:

$$\langle \mathbf{e}^1, \mathbf{y} \rangle \geq 0, \ \langle \hat{\mathbf{y}}^{(s)}, \mathbf{y} \rangle \geq 0,$$

where $\hat{\mathbf{y}}^{(s)} = (w_2^{(s)}, w_1^{(s)})^T$. Here by $\langle \mathbf{a}, \mathbf{b} \rangle$ we denote a scalar product of two vectors $\mathbf{a}, \mathbf{b} \in R^2$, i.e. $\langle \mathbf{a}, \mathbf{b} \rangle = a_1 b_1 + a_2 b_2$. For any arbitrary nonequal vectors $\mathbf{y}, \mathbf{y}^*$ the inclusion $\mathbf{y} - \mathbf{y}^* \in M_s$ is equivalent to the system of the inequalities $\langle \mathbf{e}^1, \mathbf{y} - \mathbf{y}^* \rangle \geq 0$, $\langle \hat{\mathbf{y}}^{(s)}, \mathbf{y} - \mathbf{y}^* \rangle \geq 0$, where at least one of this inequalities is strict. From here one can conclude that the set of nondominated vectors $\text{Ndom}_{M_s}(Y)$ is the image of the Pareto-optimal solutions $P_\mathbf{g}(X)$ in the multicriteria problem with the criteria $\mathbf{g}$ under mapping $\mathbf{f}$, where $\mathbf{g}$ has the following components $g_1 = f_1$, $g_2 = w_2^{(s)} f_1 + w_1^{(s)} f_2$. It means $\hat{P}_\mathbf{g}(Y) = f(P_\mathbf{g}(X))$. Thus, the set $\hat{P}_\mathbf{g}(Y)$ is the group choice according the majority preference relation $\succ$ and given "quanta" of information. Theorem 1 is proved.

Remark that in the opposite case to Theorem 1, when all second components of the vectors $\mathbf{y}^{(1)}$, $\mathbf{y}^{(2)}$, and $\mathbf{y}^{(3)}$ are positive, and all first components are negative, one should renumber the components.

Now, consider situation (II): only two vectors have positive component with the same index.

**Theorem 2.** *Let situation (II) is hold, and the vectors*

$$\mathbf{y}^{(1)} = \begin{pmatrix} w_1^{(1)} \\ -w_2^{(1)} \end{pmatrix}, \ \mathbf{y}^{(2)} = \begin{pmatrix} w_1^{(2)} \\ -w_2^{(2)} \end{pmatrix}, \ \text{and } \mathbf{y}^{(3)} = \begin{pmatrix} -w_1^{(3)} \\ w_2^{(3)} \end{pmatrix}$$

*such that* $\mathbf{y}^{(1)} \succ_1 \mathbf{0}_2$, $\mathbf{y}^{(2)} \succ_2 \mathbf{0}_2$, *and* $\mathbf{y}^{(3)} \succ_3 \mathbf{0}_2$ *are given. It is assumed that vectors* $\mathbf{y}^{(1)}$ *and* $\mathbf{y}^{(2)}$ *are not codirectional. Then* $\hat{P}_\mathbf{g}(Y) \subseteq P(Y)$, *where* $\hat{P}_\mathbf{g}(Y) = \mathbf{f}(P_\mathbf{g}(X))$, *and* $g_1 = f_1$, $g_2 = w_2^{(1)} f_1 + w_1^{(1)} f_2$, *if* $W(1, 2) > 0$ *is valid, or* $g_2 = w_2^{(2)} f_1 + w_1^{(2)} f_2$, *if* $W(1, 2) < 0$ *is valid.*

**Proof.** Let $M_l = \text{cone}\{\mathbf{e}^1, \mathbf{e}^2, \mathbf{y}^{(l)}\} \setminus \{\mathbf{0}_2\}$ for any $l \in \{1, 2, 3\}$. According to the proof of Theorem 1 in the similar way we can obtain that if $W(1, 2) > 0$, then the inclusion $M_1 \subset M_2$ is hold that implies $M_1 \cap M_2 = M_1$. And in case of $W(1, 2) < 0$ the inclusion $M_2 \subset M_1$ is valid that implies $M_1 \cap M_2 = M_2$.

One can prove that for the cones $M_1$ and $M_3$ the intersection $M_1 \cap M_3$ is equal to nonnegative orthant $R_+^2$, also $M_2 \cap M_3 = R_+^2$. Check this statement for the cones $M_1$ and $M_3$, for the cones $M_2$ and $M_3$ it can be done in a similar way. The cone $M_1$ coincides with the set of all nonzero solutions of the following system of linear inequalities:

$$\langle \mathbf{e}^1, \mathbf{y} \rangle \geq 0, \ \langle \hat{\mathbf{y}}^{(1)}, \mathbf{y} \rangle \geq 0, \quad (5)$$

where $\hat{\mathbf{y}}^{(1)} = (w_2^{(1)}, w_1^{(1)})^T$, and the cone $M_3$ – with the system

$$\langle \mathbf{e}^2, \mathbf{y} \rangle \geq 0, \ \langle \hat{\mathbf{y}}^{(3)}, \mathbf{y} \rangle \geq 0, \quad (6)$$

where $\hat{\mathbf{y}}^{(3)} = (w_2^{(3)}, w_1^{(3)})^T$. For system of inequalities (5), (6) one can find the fundamental set of solutions, which consists of the vectors $\mathbf{e}^1, \mathbf{e}^2$. For that reason $M_1 \cap M_3 = R_+^2$.

As a result, we have $M = M_1 \cup R_+^2 \cup R_+^2 = M_1$, if $W(1, 2) > 0$, and $M = M_2 \cup R_+^2 \cup R_+^2 = M_2$ in case $W(1, 2) < 0$. It means that in situation (II) the cone $M$ is also convex. Thus, we can conclude that the relation $\succ_M$ is a transitive part of the majority preference relation $\succ$, $M \subseteq K$. Using analogous arguments as in Theorem 1 one can ends the proof of the theorem. Theorem 2 is proved.

Remark that if among the vectors $\mathbf{y}^{(1)}$, $\mathbf{y}^{(2)}$, and $\mathbf{y}^{(3)}$ there are two vectors with second positive component and first negative component, one should renumber the components to use Theorem 2.

*4.2. Case of two "quanta" of information for each DM*

Now consider the situation, when each DM has two "quanta" of information. Thus, let there exist such vectors $\mathbf{y}^{(1)} = (w_1^{(1)}, -w_2^{(1)})^T$, $\mathbf{y}^{(2)} = (w_1^{(2)}, -w_2^{(2)})^T$, and $\mathbf{y}^{(3)} = (w_1^{(3)}, -w_2^{(3)})^T$, which are not codirectional, that $\mathbf{y}^{(1)} \succ_1 \mathbf{0}_2$, $\mathbf{y}^{(2)} \succ_2 \mathbf{0}_2$, and $\mathbf{y}^{(3)} \succ_3 \mathbf{0}_2$, and such vectors $\overline{\mathbf{y}}^{(1)} = (-v_1^{(1)}, v_2^{(1)})^T$, $\overline{\mathbf{y}}^{(2)} = (-v_1^{(2)}, v_2^{(2)})^T$, and $\overline{\mathbf{y}}^{(3)} = (-v_1^{(3)}, v_2^{(3)})^T$, which are not codirectional, that $\overline{\mathbf{y}}^{(1)} \succ_1 \mathbf{0}_2$, $\overline{\mathbf{y}}^{(2)} \succ_2 \mathbf{0}_2$, and $\overline{\mathbf{y}}^{(3)} \succ_3 \mathbf{0}_2$. Let for any $l \in \{1, 2, 3\}$ the condition

$$w_1^{(l)} v_2^{(l)} - w_2^{(l)} v_1^{(l)} > 0 \qquad (7)$$

is valid. The implementation of formulae (7) means that for any $l \in \{1, 2, 3\}$ the "quanta" of information of the $DM_l$ given by the vectors $\mathbf{y}^{(1)}$, $\mathbf{y}^{(2)}$, and $\mathbf{y}^{(3)}$ are consistent[12]. Denote $V(l, s) = v_2^{(l)} v_1^{(s)} - v_1^{(l)} v_2^{(s)}$.

**Theorem 3.** *Let there be two consistent "quanta" of information of the $DM_l$ given by the vectors $\mathbf{y}^{(l)}$, $\overline{\mathbf{y}}^{(l)}$ for any $l \in \{1, 2, 3\}$, as defined above. It is existed the indices $s_1, s_2 \in \{1, 2, 3\}$ such the inequalities $W(l_1, s_1) > 0$, $W(s_1, k_1) > 0$, $V(l_2, s_2) > 0$, and $V(s_2, k_2) > 0$ are hold for some indices $l_q, s_q, k_q \in \{1, 2, 3\}$, $l_q \neq k_q$, $l_q \neq s_q$, $k_q \neq s_q$ for any $q \in \{1, 2\}$. Also assume that $w_1^{(s_1)} v_2^{(s_2)} - w_2^{(s_1)} v_1^{(s_2)} > 0$ is valid. Then $\hat{P}_{\mathbf{g}}(Y) \subseteq P(Y)$, where $\hat{P}_{\mathbf{g}}(Y) = \mathbf{f}(P_{\mathbf{g}}(X))$, and $g_1 = v_2^{(s_2)} f_1 + v_1^{(s_2)} f_2$, $g_2 = w_2^{(s_1)} f_1 + w_1^{(s_1)} f_2$. Note that there are only ones such indices $s_1$, $s_2$ for particular components of vectors $\mathbf{y}^{(1)}$, $\mathbf{y}^{(2)}$, and $\mathbf{y}^{(3)}$.*

**Proof.** Denote by $M_l$ a cone generated by the vectors $\mathbf{e}^1$, $\mathbf{e}^2$, $\mathbf{y}^{(l)}$, and $\overline{\mathbf{y}}^{(l)}$ (without the origin $\mathbf{0}_2$), i. e. $M_l = \text{cone}\{\mathbf{e}^1, \mathbf{e}^2, \mathbf{y}^{(l)}, \overline{\mathbf{y}}^{(l)}\} \setminus \{\mathbf{0}_2\}$, which is convex and pointed due to (7). The cone $M_l$ can be considered as union of two cones $M_l^1 = \text{cone}\{\mathbf{e}^1, \mathbf{e}^2, \mathbf{y}^{(l)}\} \setminus \{\mathbf{0}_2\}$ and $M_l^2 = \text{cone}\{\mathbf{e}^1, \mathbf{e}^2, \overline{\mathbf{y}}^{(l)}\} \setminus \{\mathbf{0}_2\}$. According to the proof of Theorem 1 for any $l, s \in \{1, 2, 3\}$, $l \neq s$, it is valid either inclusion $M_l^q \subset M_s^q$ or inclusion $M_s^q \subset M_l^q$, where $q = 1, 2$. From here one can conclude that for any $q \in \{1, 2\}$ the cones $M_1^q$, $M_2^q$, and $M_3^q$ satisfy the inclusions $M_{l_q}^q \subset M_{s_q}^q \subset M_{k_q}^q$, where indices $l_q, s_q, k_q \in \{1, 2, 3\}$, $l_q \neq k_q$, $l_q \neq s_q$, $k_q \neq s_q$. Evidently, in such case the inequalities $W(l_1, s_1) > 0$, $W(s_1, k_1) > 0$, $V(l_2, s_2) > 0$, and $V(s_2, k_2) > 0$ are hold.

As a result, we have that unions of intersections $M = (M_1 \cap M_2) \cup (M_1 \cap M_3) \cup (M_2 \cap M_3) = M_{s_1}^1 \cup M_{s_2}^2$. It is easy to see that cones $M_{s_1}^1$, $M_{s_2}^2$ are convex, but $M$ not should be the same. Similarly to the proof of Lemma 3 one can obtain that the cone $M$ is convex if and only if $w_1^{(s_1)} v_2^{(s_2)} - w_2^{(s_1)} v_1^{(s_2)} > 0$.

It implies that the relation $\succ_M$ is a transitive part of the majority preference relation $\succ$, $M \subseteq K$. One can obtain that the cone $M$ coincides with the set of all nonzero solutions of the following system of linear inequalities:

$$\langle \hat{\mathbf{y}}^{(1)}, \mathbf{y} \rangle \geqq 0, \; \langle \hat{\mathbf{y}}^{(2)}, \mathbf{y} \rangle \geqq 0,$$

where $\hat{\mathbf{y}}^{(1)} = (v_2^{(s_2)}, v_1^{(s_2)})^T$, $\hat{\mathbf{y}}^{(2)} = (w_2^{(s_1)}, w_1^{(s_1)})^T$.

The end of the proof is similar to the proof of Theorem 1. Theorem 3 is proved.

**Theorem 4.** *If in Theorem 3 the inequality $w_1^{(s_1)} v_2^{(s_2)} - w_2^{(s_1)} v_1^{(s_2)} > 0$ is not valid, then there exist two vectors $\mathbf{u}^{(1)} = (u_1^{(1)}, -u_2^{(1)})^T$, $\mathbf{u}^{(2)} = (-u_1^{(2)}, u_2^{(2)})^T$, $\mathbf{u}^{(1)}, \mathbf{u}^{(2)} \in K$ with components $u_1^{(1)}, u_2^{(1)}, u_1^{(2)}, u_2^{(2)} > 0$, which satisfy the inequality $u_1^{(1)} u_2^{(2)} - u_2^{(1)} u_1^{(2)} > 0$ and at least one of the conditions $u_1^{(1)} v_2^{(s_2)} - u_2^{(1)} v_1^{(s_2)} > 0$ or $w_1^{(s_1)} u_2^{(2)} - w_2^{(s_1)} u_1^{(2)} > 0$. Then the inclusion $\hat{P}_{\mathbf{g}}(Y) \subseteq P(Y)$ is hold. Here $\hat{P}_{\mathbf{g}}(Y) = \mathbf{f}(P_{\mathbf{g}}(X))$, and $g_1 = u_2^{(2)} f_1 + u_1^{(2)} f_2$, $g_2 = u_2^{(1)} f_1 + u_1^{(1)} f_2$. Note that the components of the vectors $\mathbf{u}^{(1)}$, $\mathbf{u}^{(2)}$ are not unique.*

In this case the cone $M$ is not convex, and the proof is based on Lemma 3. Due to the way we choose a convex cone in the cone $M$, the Pareto set $\hat{P}_{\mathbf{g}}(Y)$ as a solution of the group choice problem is not unique. And all such Pareto sets are equally "optimal" for the group of the DMs according to the given "quanta" of information about the preference relations $\succ_1$, $\succ_2$, $\succ_3$. This is the difference between the result of Theorem 4 and the results of Theorems 1–3.

## 5. Example

Consider the example. Let $m = 2$, $n = 3$, i. e. there are 3 DMs, and $Y \subseteq R^2$. And let each $DM_l$ gives two "quanta" of information: there are vectors $\mathbf{y}^{(l)}$, $\overline{\mathbf{y}}^{(l)}$ such $\mathbf{y}^{(l)} \succ_l \mathbf{0}_2$, $\overline{\mathbf{y}}^{(l)} \succ_l \mathbf{0}_2$ for any $l \in \{1, 2, 3\}$, where

$$\mathbf{y}^{(1)} = \begin{pmatrix} 2 \\ -1 \end{pmatrix}, \mathbf{y}^{(2)} = \begin{pmatrix} 1 \\ -2 \end{pmatrix}, \mathbf{y}^{(3)} = \begin{pmatrix} 3 \\ -1 \end{pmatrix}, \overline{\mathbf{y}}^{(1)} = \begin{pmatrix} -1 \\ 1 \end{pmatrix}, \overline{\mathbf{y}}^{(2)} = \begin{pmatrix} -1 \\ 3 \end{pmatrix}, \text{ and } \overline{\mathbf{y}}^{(3)} = \begin{pmatrix} -1 \\ 4 \end{pmatrix}.$$

One can check that the conditions of Theorem 3 are hold. According to the notation we have $l_1 = 3$, $s_1 = 1$, $k_1 = 2$, $l_2 = 3$, $s_2 = 2$, $k_2 = 1$; $W(3, 1) = 1$, $W(1, 2) = 3$, $V(3, 2) = 1$, $V(2, 1) = 2$. Also the inequality $w_1^{(1)} v_2^{(2)} - w_2^{(1)} v_1^{(2)} = 5 > 0$ is valid. Then the components of vector criteria $\mathbf{g}$ are following: $g_1 = 3f_1 + f_2$, $g_2 = f_1 + 2f_2$. Let there be the set of feasible vectors $Y$, which consists of vectors $\mathbf{u}^{(1)} = (2, 3)^T$, $\mathbf{u}^{(2)} = (1, 4)^T$, $\mathbf{u}^{(3)} = (5, 2)^T$. It easy to see that $P(Y) = Y$. According to the formula of vector criteria $\mathbf{g}$, we obtain $G = \mathbf{g}(X) = \{(9, 8)^T, (7, 9)^T, (17, 9)^T\}$. From here, the Pareto set $\hat{P}_{\mathbf{g}}(Y) = \mathbf{f}(P_{\mathbf{g}}(X))$ consists of vector $\mathbf{u}^{(3)}$, which forms a group choice according to given "quanta" of information.

## 6. Conclusion

A model of multicriteria group choice problem is considered in the paper. It is assumed that each DM's preference relation is a cone relation with convex pointed cone, which contains nonnegative orthant, and does not contain the origin. Individual preferences are given in terms of "quanta" of information, reflecting the compromise, which could be done between two criteria or two groups of criteria. To aggregate individual preferences it is used the majority preference relation based on the DMs' preference relations, and its properties are investigated. It is proved that the cone of the majority preference relation, in general, is not convex, and the relation is not transitive. In case of three DMs and two criteria it is shown how to use "quanta" of information of each DM for the group choice. The set of nondominated vectors according to the majority preference relation is a Pareto set in multicriteria problem with "new" vector criteria. In one particular situation when the cone of this relation is not convex the Pareto set of "new" problem is not unique.

**Acknowledgements**

This work was supported by the Russian Foundation for Basic Research, project no. 14–07–00899.

**References**

1. Arrow, K. J. *Social Choice and Individual Values*. Wiley, New York, 1963.
2. Aizerman, M.A., Aleskerov, F.T. Vybor variantov: osnovy teorii (Variants choice: basis of the theory). Moscow: Nauka, 1990, in Russian.
3. Aleskerov, F.T., Habina, E.L., Shvarts, D.A. *Binarnie otnoshenia, grafy i kollektivnye rechenia (Binary relations, graphs and collective decisions)*. Moscow: Izd. dom GU Vyschai Shkola Ekonomiki, 2006, in Russian.
4. Mirkin, B.G. *Problema gruppovogo vybora (Problem of group choice)*. Moscow: Izdatel'stvo "Nauka", 1974, in Russian.
5. Petrovskii, A.B. *Teoriya prinyatiya reshenii (Theory of Decision Making)*. Moscow: Akademia, 2009, in Russian.
6. Hwang Ch.-L., Lin M.-J. *Group decision making under multiple criteria. Methods and applicatoins*. Springer-Verlag, 1987.
7. Kilgour, D.M., Eden, C. *Handbook of Group Decision and Negotiation*. Springer, Dordrecht, 2010.
8. Ehrgott, M., Figueira, J.R., Greco S. *Trends in Multiple Criteria Decision Analysis*. Springer, 2010.
9. Noghin, V.D. *Prinatie reshenii v mnogokriterial'noi srede: collichestvennyi podkhod (Decision Making in Multicriteria Sphere: Quantitative Approach)*. Moscow: Fizmatlit, 2005, in Russian.
10. Noghin, V.D. An algorithm of Pareto set reducing based on arbitrary finite collection of «quanta» of information. *Iskusstvennyi intellekt i prinyatie reshenii*, 2013, No. 1, pp. 63–69, in Russian.
11. Noghin V.D. Axiomatic approach to reduce the Pareto set: computational aspects. Moscow international conference on operational research (ORM2013), pp. 58–60.

12. Klimova, O.N., Noghin V.D. Using interdependent information on the relative importance of criteria in decision making. C*omputational Mathematics and Mathematical Physics,* 2006, Vol. 46, No. 12, pp. 2080–2091.